\documentstyle[12pt]{article}

\author{Gustav W. Delius
\\
\small{Department of Mathematics, King's College London}\\
\small{Strand, London WC2R 2LS, UK}\\
\footnotesize{e-mail: delius@mth.kcl.ac.uk}
\\ \\
Mark D. Gould
\\
\small{Department of Mathematics, University of Queensland}\\
\small{Brisbane Qld 4072, Australia.}\\ \\
Yao-Zhong Zhang
\\
\small{Yukawa Institute for Theoretical Physics}\\
\small{Kyoto University, Kyoto 606, Japan}\\
\footnotesize{e-mail: yzzhang@yukawa.kyoto-u.ac.jp} }

\title{Twisted Quantum Affine Algebras\\
and Solutions to the Yang-Baxter Equation}

\makeatletter

\newfont{\Bbb}{msbm10 scaled 1\@ptsize00}
\newcommand{\CC}{\mbox{\Bbb C}}
\newcommand{\ZZ}{\mbox{\Bbb Z}}

\newif\if@fewtab\@fewtabtrue

\def\draftdate{\number\day.\number\month.\number\year\ \ \ \hourmin }
{\count255=\time\divide\count255 by 60
\xdef\hourmin{\number\count255}
\multiply\count255 by-60\advance\count255 by\time
\xdef\hourmin{\hourmin:\ifnum\count255<10 0\fi\the\count255}}
\def\ps@draft{\let\@mkboth\@gobbletwo
    \def\@oddhead{}
    \def\@oddfoot
       {\hbox to 7 cm{$\scriptstyle\bf Draft\ version:\ \draftdate$
       \hfil}\hskip -7cm\hfil\rm\thepage \hfil}
    \def\@evenhead{}\let\@evenfoot\@oddfoot}

\def\label#1{\ifnum\draftcontrol=1
 \global\def\draftnote{\scriptsize\tt #1}\fi
 \@bsphack\if@filesw {\let\thepage\relax
   \def\protect{\noexpand\noexpand\noexpand}%
\xdef\@gtempa{\write\@auxout{\string
      \newlabel{#1}{{\@currentlabel}{\thepage}}}}}\@gtempa
   \if@nobreak \ifvmode\nobreak\fi\fi\fi
  \@esphack}

\def\@eqnnum{\hbox to 3cm{\phantom{\rm(\theequation)} \draftnote
                         \hfil}\hskip -3cm {\rm(\theequation)}}

\def\eqnarray{\def\draftnote{{}}\global\@fewtabtrue
\stepcounter{equation}\let\@currentlabel=\theequation
\global\@eqnswtrue
\global\@eqcnt\z@\tabskip\@centering\let\\=\@eqncr
$$\halign to \displaywidth\bgroup\@eqnsel\hskip\@centering\@eqcnt\z@
  $\displaystyle\tabskip\z@{##}$&\global\@eqcnt\@ne
  \hskip 1\arraycolsep \hfil${##}$\hfil
  &\global\@eqcnt\tw@ \hskip 1\arraycolsep
$\displaystyle\tabskip\z@{##}$
\hfil  \tabskip\@centering&\global\@eqcnt\thr@@\llap{##}\tabskip\z@
\cr}

\def\endeqnarray{\@@eqncr\egroup
      \global\advance\c@equation\m@ne$$\global\@ignoretrue}

\def\@@eqncr{\let\@tempa\relax
    \ifcase\@eqcnt \def\@tempa{& & &}\or \def\@tempa{& &}
      \or \def\@tempa{&}
      \or\def\@tempa{}
\fi\@tempa
\if@eqnsw
\if@fewtab\@eqnnum\fi
\stepcounter{equation}\fi\global
\@eqnswtrue\global\@eqcnt\z@\global\@fewtabtrue\cr}

\@addtoreset{equation}{section}

\def\cases#1{\left\{\,\vcenter{\normalbaselines\m@th
    \ialign{$\displaystyle{##}\hfil$&\quad##\hfil\crcr#1\crcr}}\right.}

\def\ct#1{\ifnum\draftcontrol=1{\tt [#1]}\else{\cite{#1}}\fi}
\def\ctz#1#2{\ifnum\draftcontrol=1{\tt [#1,#2]}\else{\cite[#1]{#2}}\fi}

\def\draftcite#1{\ifnum\draftcontrol=1#1\else{}\fi}

\def\@lbibitem[#1]#2{\item{}\hskip -3cm \hbox to 2cm
{\hfil$\scriptstyle\draftcite{#2}$}\hskip
1cm[\@biblabel{#1}]\if@filesw
     {\def\protect##1{\string ##1\space}\immediate
      \write\@auxout{\string\bibcite{#2}{#1}}}\fi\ignorespaces}

\def\@bibitem#1{\item\hskip -3cm \hbox to 2cm
{\hfil \scriptsize\tt\draftcite{#1}}\hskip 1cm
\if@filesw \immediate\write\@auxout
       {\string\bibcite{#1}{\the\value{\@listctr}}}\fi\ignorespaces}

\makeatother

\def\lb#1{\label{#1}}
\def\lab#1{\ifnum\draftcontrol=1{{\tt [#1]} \lb{#1}}\else{\lb{#1}}\fi}
\def\Eq#1{(\ref{#1})}
\def\theequation{{\thesection.\arabic{equation}}}
\def\[{\begin{eqnarray}}

\def\non{\nonumber \\ }
\def\]{\end{eqnarray}}

\def\een{\end{enumerate}}
\def\ben{\begin{enumerate}}

\renewcommand{\a}{\alpha}

\renewcommand{\d}{\delta}
\newcommand{\D}{\Delta}
\newcommand{\e}{\epsilon}
\renewcommand{\e}{\epsilon}
\renewcommand{\l}{\lambda}
\renewcommand{\L}{\Lambda}
\newcommand{\s}{\sigma}

\newcommand{\half}{\frac{1}{2}}
\newcommand{\<}{\langle}
\renewcommand{\>}{\rangle}

\def\binco#1#2{\left(\begin{array}{c}#1\\#2\end{array}\right)}
\newtheorem{theorem}{Theorem}

\newcommand{\lie}{L}
\newcommand{\tlie}{\hat{\lie}^{(2)}}
\newcommand{\rootb}{\bar{\a}}
\newcommand{\ev}{\mbox{ev}}
\newcommand{\tp}{{\bf \check{P}}}
\newcommand{\ttpg}{\tilde{G}}
\newcommand{\ettpg}{\tilde{\Gamma}}
\newcommand{\imp}{\subseteq}
\newcommand{\minimal}{$\lie_0$-irreducible }

\begin{document}

\def\draft{\pagestyle{draft}\thispagestyle{draft}
\global\def\draftcontrol{1}}
\global\def\draftcontrol{0}


\begin{titlepage}
\pagestyle{empty}

\maketitle
\begin{abstract}

We construct spectral parameter dependent R-matrices for
the quantized enveloping algebras of twisted affine
Lie algebras. These give new solutions to the
spectral parameter dependent quantum Yang-Baxter equation.

\end{abstract}
\vspace{1cm}
q-alg/9508012\\
KCL-TH-95-8\\
YITP/K-1119
\end{titlepage}
\newpage


\section{Introduction\lab{sectintro}}

The solutions to the Yang-Baxter equation play a central role in the
theory of quantum integrable models \ct{Skl79,Kor93}.
In statistical mechanics they are the
Boltzman weights of exactly solvable lattice models \ct{Bax82}. In
quantum field theory they give the exact factorizable scattering
matrices \ct{Zam79}. For an introduction to the mathematical aspects
of the Yang-Baxter equation see e.g. \ct{Jim89}.

The Yang-Baxter equation with spectral parameter has the form
\[\label{YB}
R^{12}(u)R^{13}(uv)R^{23}(v)
  =R^{23}(v)R^{13}(uv)R^{12}(u).
\]
The $R^{ab}(u)$ are matrices which depend on a spectral parameter
$u$ and which act on the tensor product
of two vector spaces $V_a$ and $V_b$
\[
R^{ab}(u):~V_a\otimes V_b\rightarrow V_a\otimes V_b.
\]
The products of $R$'s in \Eq{YB} act on the space $V_1\otimes V_2
\otimes V_3$.

The mathematical framework for the construction of trigonometric
solutions of the quantum Yang-Baxter equation \Eq{YB}
is given by the quantum affine algebras $U_q(\hat{\lie})$ introduced
by Jimbo \ct{Jim85}
and  Drinfeld \ct{Dri85}. These are deformations of the enveloping
algebras $U(\hat{\lie})$ of affine Lie algebras \ct{Kac90}.
Associated to any two finite-dimensional irreducible
$U_q(\hat{\lie})$-modules $V(\l)$ and $V(\mu)$ there exists a trigonometric
R-matrix $R^{\l\mu}(u)$. Given three modules, the R-matrices for all
pairs of these three modules are a solution of \Eq{YB}.
Many R-matrices of untwisted quantum affine algebras have since been determined
(see references in \ct{Del94a}), leading to a large number of new quantum
integrable models, quantum spin chains, exactly solvable lattice models and
exact scattering
matrices. The method has also been extended to quantum affine superalgebras
\ct{Del95a}.

While it was clear from the beginning \ct{Dri85} that all (twisted and
untwisted) affine Lie algebras can be quantized and give rise to
trigonometric R-matrices, the twisted algebras have hardly been treated
in the literature. The only R-matrices associated to twisted quantum
affine algebras which we have found in the literature are those associated
to the vector representation of $U_q(A_l^{(2)})$ and $U_q(D_{l+1}^{(2)})$
\ct{Baz85,Jim86}. (The $U_q(A_2^{(2)})$ R-matrix was found already in
\ct{Ize81}; another R-matrix for $U_q(D^{(2)}_{l+1})$ has been found
in \cite{Grimm} by ``Baxterizing"  the so-called dilute BWM algebra
\footnote{We thank Ole Warnaar for drawing our attention to the ref.
\cite{Grimm}}.)
The knowledge of these R-matrices has had many physical
applications. They have for example been used to obtain transfer
matrices of solvable lattice models \ct{Kun94} or to
diagonalize quantum spin chain Hamiltonians on the periodic chain
\ct{Res87} and on the open chain \ct{Art94,Art95}.
To generalize these works to the models associated to higher dimensional
representations it is necessary to know the corresponding R-matrices and
that is the topic of this paper.

R-matrices are also needed to construct the S-matrices of quantum field
theories with quantum affine symmetries \ct{Del95c}. In particular
the R-matrices associated to all fundamental representations of
$A_{2n}^{(2)}$ which we construct in this paper will give the
scattering matrices for the solitons in $A_{2n}^{(2)}$ affine Toda
theory.

The paper is organized as follows:
In section \ref{secttakm} we review the necessary facts about twisted
Lie affine algebras \ct{Kac90}. In section \ref{secttqaa} we discuss the
quantized affine algebras and give the equations which uniquely determine
the R-matrices (Jimbo's equations). In section \ref{sectjimbo} we
explain how to solve these equations. Our technique is an
extension of the tensor product
graph method introduced in \ct{Zha91} and generalized in \ct{Del94a}.
We have obtained the R-matrices for the following twisted algebras
and tensor products:
\[
\begin{array}{|cc|cc|l|}
\hline
U_q(A_{2l}^{(2)}) & l\geq 2 & V(\l_k)\otimes V(\l_r) & k+r\leq l &
\mbox{section \ref{sectjimbo1}} \\
U_q(A_{2l-1}^{(2)}) & l\geq 3 & V(k\l_1)\otimes V(r\l_1) &  &
\mbox{section \ref{sectjimbo2}} \\
U_q(D_{l+1}^{(2)}) & l\geq 2 & V(k\l_l)\otimes V(r\l_l) &  &
\mbox{section \ref{sectjimbo3}} \\
\hline
\end{array}
\]
Here $\l_i$ denotes the i-th fundamental weight.
We give some of the technical
details in appendices. In particular in appendix
\ref{app2} we derive the tensor product decompositions and branching
rules which we need in the paper.


\section{Twisted affine Lie algebras\lab{secttakm}}

We recall the relevant information about twisted affine
Lie algebras \ct{Kac90}.
Let $\lie$ be a finite dimensional simple Lie algebra and
$\s$ a diagram automorphism
of $\lie$ of order $k$. Associated to these one constructs the
twisted affine Lie algebra $\hat{L}^{(k)}$. In this paper we will
assume\footnote{The only diagram automorphism which is not of order
$k=2$ is the triality of the Lie algebra $D_4$ and we will not
treat this case in this paper.} $k=2$.
Let $\lie_0$ be
the fixed point subalgebra under the diagram automorphism $\s$.
We recall that
\[\label{liesplit}
\lie=\lie_0\oplus\lie_1,~~~~
[\lie_i,\lie_j]=\lie_{(i+j) \rm{mod}2}.
\]
$\lie_1$ gives rise to an irreducible $\lie_0$-module under
the adjoint action of $\lie_0$.
Let $\theta_0$ be its highest
weight. In table \ref{tab1} we list all the cases with $k=2$.
Below we restrict ourselves to the three families and leave
out the exceptional case for technical reasons which will become
apparent later.

\begin{table}[htb]
\[
\begin{array}{cc|c|c}
\lie&&\lie_0&\theta_0\\
\hline
A_{2l},&l\geq 1 & B_l & 2\l_1=2\e_1 \\
A_{2l-1},&l\geq 3 & C_l & \l_2=\e_1+\e_2 \\
D_{l+1},&l\geq 2 & B_l & \l_1=\e_1 \\
E_6 && F_4 & \l_4
\end{array}\nonumber
\]
\caption{Table of the finite dimensional simple Lie algebras
$\lie$ which posess a diagram automorphism of order $k=2$,
their fixed point subalgebras $\lie_0$ and the highest weight
$\theta_0$ of the adjoint $\lie_0$-module $\lie_1$.
Here and in the rest of the paper we give weights
either as integer combinations of the fundamental weights
$\l_i$ or alternatively we give them in terms of the $\e_i$ which
form a basis of the root space of $gl(n)$ into which we imbed
the other algebras, see Appendix A.
\label{tab1}}
\end{table}

We recall that $\lie$ admits generators $E_i,F_i,H_i,~0\leq i\leq l,$
satisfying the defining relations
\footnote{The rescaled generators $E_i^\prime=\sqrt{2/\rootb_i^2}E_i,~
F_i^\prime=\sqrt{2/\rootb_i^2}E_i,~
H_i^\prime=2/\rootb_i H_i$ satisfy the more
usual commutation relations with the structure constants given by
the Cartan matrix.}
\[\label{lierels}
&&
[H_i,E_j]=(\rootb_i,\rootb_j)E_j,~~~
[H_i,F_j]=-(\rootb_i,\rootb_j)F_j,~~~
[E_i,F_j]=\delta_{ij}H_i,
\non &&
(\mbox{ad}_{E_i})^{1-a_{ij}}E_j=
(\mbox{ad}_{F_i})^{1-a_{ij}}F_j=0,~~~i\neq j,
\]
where $a_{ij}=2(\rootb_i,\rootb_j)/(\rootb_i,\rootb_i)$ are the entries
of the corresponding (twisted) Cartan matrix of $\tlie$.
Here the $E_i,F_i,H_i,~1\leq i\leq l,$ form the Chevalley generators
for $\lie_0$ and the $\rootb_i,~(1\leq i\leq l)$ are the simple
roots of $\lie_0$.
$E_0\in\lie_1$ corresponds to the minimal weight vector and thus has
weight $-\theta_0$. It follows that $\rootb_0=-\theta_0$ and that
$H_0=-\sum_{i=1}^l a_i H_i$ lies in the
Cartan subalgebra $H$ of $\lie_0$. The integers $a_i$ are known as
the Kac labels of $\tlie$.

Throughout we let $(~,~)$ be a fixed invariant bilinear form on
$\lie$ which induces a corresponding invariant form $(~,~)$ on
$H^*$. A suitable choice for the invariant form on $\lie$ togeher
with a realization of the simple generators is given in
Appendix \ref{app1} for completeness. With our choice we have
\[\label{form}
(E_i,F_j)=\delta_{ij},~~~(H_i,H_j)=(\rootb_i,\rootb_j).
\]

We now introduce the corresponding twisted affine Lie algebra
$\tlie{}^\prime$ which admits the decomposition
\[
\tlie{}^\prime=\bigoplus_{m\in\half\ZZ}\hat{\lie}_m
\oplus \CC c_0,~~~~
\hat{\lie}_m=\left\{
\begin{array}{l}
\lie_0(m),~m\in\ZZ\\
\lie_1(m),~m\in\ZZ+\half
\end{array}\right.
\]
with $\lie_a(m)=\{x(m)|x\in\lie_a\},~a=0,1$
and $c_0$ a central charge.
The Lie bracket is given by
\[
[x(m),y(n)]=[x,y](m+n)+m~c_0~\delta_{m+n,0}~(x,y),~~~
[c_0,x(m)]=0.
\]
Here $(~,~)$ is the fixed invariant
bilinear form on $\lie$. Note that $\hat{\lie}_0=\lie_0$.
A suitable set of generators for $\tlie{}^\prime$ is given by
\[\label{simpgen}
&&e_i=E_i(0),~~~h_i=H_i(0),~~~f_i=F_i(0),~~~1\leq i\leq l,
\non
&&e_0=E_0(1/2),~~~h_0=H_0(0)+1/2 c_0,~~~f_0=F_0(-1/2).
\]
This algebra is extended to $\tlie=\tlie{}^\prime\oplus\CC d_0$
by the introduction of the level operator $d_0$ satisfying
\[
[d_0,x(m)]=m\, x(m),~~~[d_0,c_0]=0.
\]
As a Cartan subalgebra of $\tlie$ we take
\[
\hat{H}=H(0)\oplus \CC c_0\oplus \CC d_0.
\]
The weights for $\tlie$ are of the form $\l=(\bar{\l},c_\l,d_\l)$
where $\bar{\l}\in H^*$ and $c_\l, d_\l$ are the eigenvalues of
the central extension $c_0$ and the level operator $d_0$
respectively. The simple roots corresponding to the set of
simple generators in \Eq{simpgen} are
\[
\a_i=(\rootb_i,0,0),~~1\leq i\leq l,~~~
\a_0=(-\theta_0,0,\half).
\]
The invariant bilinear form on $\hat{H}^*$ is given by
\[
(\l,\mu)=(\bar{\l},\bar{\mu})+c_\l d_\mu+ d_\l c_\mu.
\]
With this convention we have $(\a_i,\a_j)=(\rootb_i,\rootb_j),~
(0\leq i,j\leq l)$ and our simple generators satisfy the
defining relations
\[
&&
[h_i,e_j]=(\a_i,\a_j)e_j,~~~
[h_i,f_j]=-(\a_i,\a_j)f_j,~~~
[e_i,f_j]=\delta_{ij}h_i,
\non &&
(\mbox{ad}_{e_i})^{1-a_{ij}}e_j=
(\mbox{ad}_{f_i})^{1-a_{ij}}f_j=0,~~~i\neq j.
\]
Following Kac \ct{Kac90} it is useful to introduce the
weights $\gamma=(\bar{0},1,0)$ and $\delta=(\bar{0},0,1)$ so
that $(\gamma,\delta)=1, (\gamma,\gamma)=(\delta,\delta)=0$.
Our simple roots are then given by $\a_i=\rootb_i,~1\leq i\leq l,~~
\a_0=-\theta_0+1/2 \delta$.


We have an algebra homomorphism, called the
{\em evaluation map}, $\ev_t:U(\tlie)\rightarrow\CC[t,t^{-1}]
\otimes U(\lie)$, with $U(\tlie),U(\lie)$ the enveloping algebras
of $\tlie, \lie$ respectively, given by
\[
\ev_t(x(m))=t^{2m}x,~~~
\ev_t(c_0)=0,~~~
\ev_t(d_0)=\half\,t\frac{d}{dt},
\]
and extended to all of $U(\tlie)$ in the natural way.
Thus given a finite dimensional $\lie$-module $V$
carrying a representation $\pi$ we have a
corresponding $\tlie$ module $V(t)=\CC[t,t^{-1}]\otimes V$
carrying the {\em loop representation} $\hat{\pi}$ given by
\[
\hat{\pi}=(1\otimes\pi)\ev_t.
\]
Below we consider the problem of quantizing such representations
to give solutions of the Yang-Baxter equation.



An important role will be played below by those irreducible $\lie$-modules
which are also irreducible under the $\lie_0$ subalgebra. We call these
{\em \minimal}
modules. We will see below that the loop representations
built on \minimal modules can all be quantized.
Such \minimal modules appear to exist for the first three
cases in table \ref{tab1} only and this is the reason why we are
restricting
to these cases. In table \ref{tab2} we list for each of the
three families the
highest weights of all the \minimal irreps of $\lie$ together with
their highest weight with respect to $\lie_0$.

\begin{table}[htb]
\[
\begin{array}{c|c|c|cc}
\lie & \lie_0 & \L & \L_0 &\\
\hline
A_{2l} & B_l & \l_k,\l_{2l+1-k} &
\l_k, & 1\leq k\leq l \\
A_{2l-1} & C_l & a \l_1 & a \l_1, & a\in \ZZ_+ \\
D_{l+1} & B_l & a \l_l,a\l_{l+1} &
a \l_l, & a\in\ZZ_+
\end{array}\nonumber
\]
\caption{\minimal irreps.
$\L$ are the highest weights of the \minimal irreps of $\lie$
and $\L_0$ are the corresponding highest weights under $\lie_0$.
\lab{tab2}}
\end{table}

In appendix \ref{app2} we show that the tensor product of any two
such \minimal $\lie$-modules decomposes into a {\em multiplicity free}
direct sum of irreducible $\lie_0$-modules. This is important
because it implies that a solution to Jimbo's equations
will always exist for such tensor products (see below).


\section{Twisted quantum affine algebras\lab{secttqaa}}

Corresponding to the twisted affine algebra $\tlie$
we have the twisted quantum affine algebra $U_q(\tlie)$ with
generators $q^{\pm h_i/2},e_i,f_i,d_0,~(0\leq i\leq l)$ and
defining relations
\[\label{qrels}
&&
[h_i,e_j]=(\a_i,\a_j)e_j,~~~
[h_i,f_j]=-(\a_i,\a_j)f_j,~~~
[h_i,h_j]=0,
\non &&
[d_0,e_i]=\half\d_{i,0}e_i,~~~
[d_0,f_i]=-\half\d_{i,0}f_i,~~~
[d_0,h_i]=0,
\non &&
[e_i,f_j]=\delta_{ij}\frac{q^{h_i}-q^{-h_i}}{q-q^{-1}},
\non &&
\sum_{k=0}^{1-a_{ij}} (-1)^k e_i^{(1-a_{ij}-k)} e_j e_i^{(k)}=0,
{}~~i\neq j,
\non &&
\sum_{k=0}^{1-a_{ij}} (-1)^k f_i^{(1-a_{ij}-k)} f_j f_i^{(k)}=0,
{}~~i\neq j,
\]
where
\[
&&
e_i^{(k)}=\frac{e_i^k}{[k]_{q_i}!},~~
f_i^{(k)}=\frac{f_i^k}{[k]_{q_i}!},
\non &&
[k]_q=\frac{q^k-q^{-k}}{q-q^{-1}},~~~
[k]_q!=\prod_{n=1}^k[n]_q,~~~
q_i=q^{\half(\a_i,\a_i)}.
\]
$U_q(\tlie)$ is a quasi-triangular Hopf algebra with coproduct
$\D$ and antipode $S$ given by
\[\label{Hopf}
\D(e_i)=q^{-h_i/2}\otimes e_i+ e_i\otimes q^{h_i/2},~~~
&&S(e_i)=-q_i e_i,
\non
\D(f_i)=q^{-h_i/2}\otimes f_i+ f_i\otimes q^{h_i/2},
&&S(f_i)=-q_i^{-1} f_i,
\non
\D(q^{\pm h_i/2})=q^{\pm h_i/2}\otimes q^{\pm h_i/2},
&&S(q^{\pm h_i/2})=q^{\mp h_i/2},
\non
\D(d_0)=1\otimes d_0+d_0\otimes 1,
&& S(d_0)=-d_0.
\]

Throughout $\bar{R}$ denotes the universal R-matrix of
$U_q(\tlie)$ which by definition satisfies
\[\label{rdef}
&&
\bar{R}\D(a)=\D^T(a)\bar{R},~~\forall a\in U_q(\tlie),
\non &&
(1\otimes\D)\bar{R}=\bar{R}_{13}\bar{R}_{12},~~~
(\D\otimes1)\bar{R}=\bar{R}_{13}\bar{R}_{23}
\]
where $\D^T(a)$ is the opposite coproduct. A direct consequence
of the above relations is that $\bar{R}$ satisfies the
quantum Yang-Baxter equation
\[\label{qYB}
\bar{R}_{12}\bar{R}_{13}\bar{R}_{23}=
\bar{R}_{23}\bar{R}_{13}\bar{R}_{12}
\]
Note that the generators $q^{\pm h_i/2}, e_i, f_i,~(1\leq i\leq l)$
generate the quantum algebra $U_q(\lie_0)$ which is a
quasi-triangular Hopf subalgebra of $U_q(\tlie)$. We denote its
universal R-matrix by $R$.

We shall see below that any minimal irrep $V_0(\l)$ of $U_q(\lie_0)$
can be affinized to give rise to an irrep of $U_q(\tlie)$.
To perform such an affinization it is necassary and sufficient
to find operators $\pi_\l(e_0)$ and $\pi_\l(f_0)$ acting on
$V_0(\l)$ which satisfy the required defining relations
\Eq{qrels} of $U_q(\tlie)$.

We define an automorphism $D_t$ of $U_q(\tlie)$ by
\[
D_t(e_i)=t^{\d_{i0}}e_i,~~
D_t(f_i)=t^{-\d_{i0}}f_i,~~
D_t(h_i)=h_i.
\]
Given any two minimal irreps $\pi_\l$ and $\pi_\mu$ of $U_q(\lie_0)$
and their affinizations to irreps of $U_q(\tlie)$,
we obtain a one-parameter family of representations
$\D_{\l\mu}^u$ of $U_q(\tlie)$ on
$V_0(\l)\otimes V_0(\mu)$ defined by
\[
\D_{\l\mu}^u(a)=\pi_\l\otimes \pi_\mu\left((D_u\otimes 1)
\D(a)\right),~~~\forall a\in U_q(\tlie),
\]
where $u$ is the spectral parameter. We define the spectral
parameter dependent R-matrix
\[
R^{\l\mu}(u)=(\pi_\l\otimes\pi_\mu)\left((D_u\otimes 1)\bar{R}\right).
\]
It follows from \Eq{qYB} that this R-matrix gives a solution to
the spectral parameter dependent Yang-Baxter equation \Eq{YB}.
{}From the defining property \Eq{rdef} of the universal R-matrix
one derives the equations
\[\label{Jimbo}
R^{\l\mu}(u)\,\D_{\l\mu}^u(a)=(\D^T)_{\l\mu}^u(a)\,R^{\l\mu}(u)
\]
which, because the representations $\D_{\l\mu}^u$
are irreducible for generic $u$,
uniquely determine $R^{\l\mu}(u)$ up to a scalar function
of $u$. These are the {\em Jimbo equations} for twisted affine
algebras.

As in \ct{Del94a} we normalize $R^{\l\mu}(u)$ such that
\[
\check{R}^{\l\mu}(u)\check{R}^{\mu\l}(u^{-1})=I~~
\mbox{and}~~R(0)=\pi_\l\otimes\pi_\mu(R),
\]
where $R$ is the R-matrix of $U_q(\lie_0)$ and
$\check{R}^{\l\mu}(u)=P\,R^{\l\mu}(u)$ with
$P:V_0(\l)\otimes V_0(\mu)\rightarrow V_0(\mu)\otimes V_0(\l)$
the usual permutation operator.

In order for the equation \Eq{Jimbo} to hold for
all $a\in U_q(\tlie)$ it is sufficient that it holds for all
$a\in U_q(\lie_0)$ and in addition for the extra generator
$e_0$.
The relation for $e_0$ reads explicitly
\[
&&
R^{\l\mu}(u)
\left(u\,\pi_\l(e_0)\otimes\pi_\mu(q^{h_0/2})+
\pi_\l(q^{-h_0/2})\otimes \pi_\mu(e_0)\right)
\non&&
=\left(u\pi_\l(e_0)\otimes\pi_\mu(q^{h_0/2})+
\pi_\l(q^{-h_0/2})\otimes\pi_\mu(e_0)\right)
R^{\l\mu}(u),
\]
or equivalently
\[
&&
\check{R}^{\l\mu}(u)
\left(u\,\pi_\l(e_0)\otimes\pi_\mu(q^{h_0/2})+
\pi_\l(q^{-h_0/2})\otimes \pi_\mu(e_0)\right)
\non&&
=\left(\pi_\mu(e_0)\otimes\pi_\l(q^{h_0/2})+
u\,\pi_\mu(q^{-h_0/2})\otimes\pi_\l(e_0)\right)
\check{R}^{\l\mu}(u).
\]


\section{Solutions to Jimbo's equations\lab{sectjimbo}}

With $V_0(\l)$ and $V_0(\mu)$ denoting two minimal irreps of
$U_q(\tlie)$ we write the tensor product decomposition into
irreducible $U_q(\lie_0)$-modules as
\[\label{dec}
V_0(\l)\otimes V_0(\mu)=\bigoplus_\nu V_0(\nu)
\]
and note that there are no multiplicities in this decomposition
for the cases which we are considering (c.f. Appendix \ref{app2}).
We let $P_\nu^{\l\mu}$ be the projection operator of
$V_0(\l)\otimes V_0(\mu)$ onto $V_0(\nu)$ and set
\[
\tp^{\l\mu}_\nu=\check{R}^{\l\mu}(1)\,P^{\l\mu}_\nu=
P^{\mu\l}_\nu\,\check{R}^{\l\mu}(1).
\]
We may thus write
\[
\check{R}^{\l\mu}(u)=\sum_\nu\,\rho_\nu(u)\,\tp^{\l\mu}_\nu,~~~
\rho_\nu(1)=1.
\]
Following our previous approach \ct{Del94a}, the coefficients
$\rho_\nu(u)$ may be determined according to the recursion
relation
\[\label{rec}
\rho_\nu(u)=\frac{q^{C(\nu)/2}+\e_\nu\e_{\nu'}u\,q^{C(\nu')/2}}
{u\,q^{C(\nu)/2}+\e_\nu\e_{\nu'}\,q^{C(\nu')/2}}
\rho_{\nu'}(u),
\]
which holds for any $\nu\neq\nu'$ for which
\[\label{edge}
P^{\l\mu}_\nu\left(\pi_\l(e_0)\otimes\pi_\mu(q^{h_0/2})\right)
P^{\l\mu}_{\nu'}\neq 0.
\]
Here $C(\nu)$ is the eigenvalue of the universal Casimir element
of $\lie_0$ on $V_0(\nu)$ and $\e_\nu$ denotes the parity
of $V_0(\nu)\imp V_0(\l)\otimes V_0(\mu)$,
 (c.f. Appendix \ref{app3}).

To graphically encode the recursive relations between the
different $\rho_\nu$ we introduce the {\bf Twisted Tensor
Product Graph} $\ttpg^{\l\mu}$ associated to the tensor
product module $V_0(\l)\otimes V_0(\mu)$. The nodes of this
graph are given by the highest weights $\nu$ of the
$U_q(\lie_0)$-modules occuring in the decomposition \Eq{dec}
of the tensor product module. There is an edge between two
nodes $\nu\neq\nu'$ iff \Eq{edge} holds.

Given a tensor product module and its decomposition, it is not
in general an easy task to determine the twisted tensor product
graph because in order to determine between which nodes of the
graph relation \Eq{edge} holds requires detailed calculations.
We therefore introduce the
{\bf Extended Twisted Tensor Product Graph} $\ettpg^{\l\mu}$
which has the same set of nodes as the twisted tensor product graph (TPG)
but has an edge between two vertices $\nu\neq\nu'$ whenever
\[\label{edge1}
V_0(\nu')\imp V_0(\theta_0)\otimes V_0(\nu)
\]
and
\[\label{edge2}
\epsilon_\nu \epsilon_{\nu'}=
\left\{\begin{array}{l}
+1~~~\mbox{if }V_0(\nu)\mbox{ and }V_0(\nu')\mbox{ are in the same
irrep of }\lie\\
-1~~~\mbox{if }V_0(\nu)\mbox{ and }V_0(\nu')\mbox{ are in different
irreps of }\lie.\end{array}\right.
\]
The conditions \Eq{edge1} and \Eq{edge2} are necessary
conditions for \Eq{edge} to hold and therefore the twisted TPG
is contained in the extended twisted TPG.
To see why \Eq{edge1} is a necessary condition for
\Eq{edge} one must
realize that $e_0\otimes q^{h_0/2}$ is the lowest component of
a tensor operator corresponding to $V_0(\theta_0)$, see \ct{Zha91}
for details.
The necessity of \Eq{edge2} follows from the following fact derived in
Appendix \ref{app3}:
Two vertices $\nu\neq\nu'$ connected by an edge in the
twisted TPG (i.e., for which \Eq{edge} is satisfied)
must have the {\em same} parity if $V_0(\nu)$ and $V_0(\nu')$
belong to the {\em same} irreducible $\lie$-module while they
must have {\em opposite} parities if they belong to different
irreducible $\lie$-modules.

While the extended twisted TPG will always include
the twisted TPG, it will in general have more
edges. Only if the extended twisted TPG is a tree are we guaranteed that it
coincides with the twisted TPG.

\noindent{\bf Note:} Unlike the untwisted case \ct{Del94a},
we may now get an edge between $\nu$ and $\nu'$ of the {\em same}
parity.
This gives rise to a twisted TPG which may
be topologically quite different to the untwisted TPG.

We will impose a relation \Eq{rec} for every edge in the
extended twisted TPG. Because the extended
TPG will in general have more edges than the
unextended twisted TPG, we will be imposing
too many relations. These relations may be inconsistent and we
are therefore not guaranteed a solution.
If however a solution exists, then it must be the
unique correct solution to Jimbo's equations.

As seen below, for the minimal cases we are considering, the
extended twisted TPG is always consistent and
thus will always give rise to a solution of the quantum Yang-Baxter
equation.

Throughout we adopt the convenient notation
\[
\<a\>_\pm=\frac{1\pm x\,q^a}{x\pm q^a},
\]
so that the relation \Eq{rec} may be expressed as
\[
\rho_\nu(u)=\left\langle\frac{C(\nu')-C(\nu)}{2}\right\>
_{\epsilon_\nu\epsilon_{\nu'}}\,
\rho_{\nu'}(u).
\]

We will now determine the R-matrices for any tensor product
of any two \minimal representations for all the three families
of twisted quantum affine algebras $U_q(A_{2l}^{(2)}),
U_q(A_{2l-1}^{(2)})$ and $U_q(D_{l+1}^{(2)})$.


\subsection{R-matrices for $U_q(A_{2l}^{(2)})$
\label{sectjimbo1}}

This is the case of the first line in table \ref{tab1}, i.e.
$\lie=A_{2l}=sl(2l+1)$, $\lie_0=B_l=
so(2l+1)$ and $\theta_0=2\e_1=2\l_1$.

The defining (vector) irrep $V(\l_1)$ of $U_q(\lie)$ is undeformed.
By this we mean that the representation matrices for the
fundamental generators are the same as those in the classical
case, i.e. they are independent of $q$. It is also a minimal
irrep, i.e. $V(\l_1)=V_0(\l_1)$ is also irreducible as a module
of $U_q(\lie_0)$. Furthermore it is affinizable, i.e. it carries
a representation of $U_q(A_{2l}^{(2)})$. Also this affinized
representation is undeformed, i.e. $\pi(e_0)$ and $\pi(f_0)$
are given by the classical expressions.

We have the corresponding twisted TPG for
$V(\l_1)\otimes V(\l_1)$
\[\label{ttpg1}
\unitlength=1mm
\linethickness{0.4pt}
\begin{picture}(102.60,12.60)(40,13)
\put(60.00,20.00){\circle*{5.20}}
\put(80.00,20.00){\circle*{5.20}}
\put(60.00,16.00){\makebox(0,0)[ct]{$0$}}
\put(80.00,16.00){\makebox(0,0)[ct]{$2\l_1$}}
\put(100.00,20.00){\circle*{5.20}}
\put(100.00,16.00){\makebox(0,0)[ct]{$\l_2$}}
\put(60.00,20.00){\line(1,0){40.00}}
\put(60.00,24.00){\makebox(0,0)[cb]{$+$}}
\put(80.00,24.00){\makebox(0,0)[cb]{$+$}}
\put(100.00,24.00){\makebox(0,0)[cb]{$-$}}
\end{picture}
\]
where $\pm$ indicate the parities. This is quite different to the
untwisted TPG
\[\label{tpg1}
\unitlength=1mm
\linethickness{0.4pt}
\begin{picture}(102.60,12.60)(40,13)
\put(60.00,20.00){\circle*{5.20}}
\put(80.00,20.00){\circle*{5.20}}
\put(60.00,16.00){\makebox(0,0)[ct]{$0$}}
\put(80.00,16.00){\makebox(0,0)[ct]{$\l_2$}}
\put(100.00,20.00){\circle*{5.20}}
\put(100.00,16.00){\makebox(0,0)[ct]{$2\l_1$}}
\put(60.00,20.00){\line(1,0){40.00}}
\put(60.00,24.00){\makebox(0,0)[cb]{$+$}}
\put(80.00,24.00){\makebox(0,0)[cb]{$-$}}
\put(100.00,24.00){\makebox(0,0)[cb]{$+$}}
\end{picture}
\]

Since $\l_2$ is an extremal node on the twisted TPG it
follows that $V_0(\l_2)$ is affinizable, i.e. it too
carries an irrep of $U_q(A_{2l}^{(2)})$ (for a discussion of
the relation between TPGs and finite dimensional
irreps of quantum affine algebras see \ct{Del95b}). More generally
we have the following twisted TPG for for
$V(\l_1)\otimes V(\l_k)$, $k<l$
\[\label{ttpg2}
\unitlength=1mm
\linethickness{0.4pt}
\begin{picture}(102.60,12.60)(40,13)
\put(60.00,20.00){\circle*{5.20}}
\put(80.00,20.00){\circle*{5.20}}
\put(100.00,20.00){\circle*{5.20}}
\put(60.00,20.00){\line(1,0){40.00}}
\put(60.00,16.00){\makebox(0,0)[ct]{$\l_{k-1}$}}
\put(80.00,16.00){\makebox(0,0)[ct]{$\l_1+\l_k$}}
\put(100.00,16.00){\makebox(0,0)[ct]{$\l_{k+1}$}}
\put(60.00,24.00){\makebox(0,0)[cb]{$+$}}
\put(80.00,24.00){\makebox(0,0)[cb]{$+$}}
\put(100.00,24.00){\makebox(0,0)[cb]{$-$}}
\end{picture}
\]
so that $\l_{k+1}$ is an extremal node and hence, by recursion,
each of the fundamental irreps $V_0(\l_k),~1\leq k\leq l$, is
affinizable. Again the above twisted TPG \Eq{ttpg2} is different
to the untwisted one which is
\[\label{tpg2}
\unitlength=1mm
\linethickness{0.4pt}
\begin{picture}(102.60,12.60)(40,13)
\put(60.00,20.00){\circle*{5.20}}
\put(80.00,20.00){\circle*{5.20}}
\put(100.00,20.00){\circle*{5.20}}
\put(60.00,20.00){\line(1,0){40.00}}
\put(60.00,16.00){\makebox(0,0)[ct]{$\l_{k-1}$}}
\put(80.00,16.00){\makebox(0,0)[ct]{$\l_{k+1}$}}
\put(100.00,16.00){\makebox(0,0)[ct]{$\l_1+\l_k$}}
\put(60.00,24.00){\makebox(0,0)[cb]{$+$}}
\put(80.00,24.00){\makebox(0,0)[cb]{$-$}}
\put(100.00,24.00){\makebox(0,0)[cb]{$+$}}
\end{picture}
\]

\noindent{\bf Note:}
In the untwisted case $V_0(\l_k),~k>1$ does not occur on extremal
nodes of any TPG and can therefore not be shown to
be affinizable to a representation of $U_q(\hat{\lie}_0^{(1)})$.
In fact it is generally not affinizable \ct{Del95b,Del94a}
in the untwisted sense. But, as seen above,
it is nevertheless affinizable to a representation of the twisted
algebra $U_q(\tlie)$.

\begin{figure}[ht]
\unitlength=1.30mm
\linethickness{0.4pt}
\begin{picture}(132.00,70.00)(17,3)
\put(30.00,10.00){\circle*{4.00}}
\put(46.00,10.00){\circle*{4.00}}
\put(38.00,18.00){\circle*{4.00}}
\put(54.00,18.00){\circle*{4.00}}
\put(130.00,10.00){\circle*{4.00}}
\put(114.00,10.00){\circle*{4.00}}
\put(122.00,18.00){\circle*{4.00}}
\put(62.00,10.00){\circle*{4.00}}
\put(46.00,26.00){\circle*{4.00}}
\put(80.00,60.00){\circle*{4.00}}
\put(72.00,52.00){\circle*{4.00}}
\put(88.00,52.00){\circle*{4.00}}
\put(96.00,44.00){\circle*{4.00}}
\put(80.00,44.00){\circle*{4.00}}
\put(30.00,10.00){\line(1,1){20.00}}
\put(46.00,10.00){\line(1,1){12.00}}
\put(46.00,26.00){\line(1,-1){16.00}}
\put(38.00,18.00){\line(1,-1){8.00}}
\put(80.00,60.00){\line(-1,-1){12.00}}
\put(80.00,60.00){\line(1,-1){20.00}}
\put(72.00,52.00){\line(1,-1){12.00}}
\put(88.00,52.00){\line(-1,-1){12.00}}
\put(130.00,10.00){\line(-1,1){12.00}}
\put(122.00,18.00){\line(-1,-1){8.00}}
\put(114.00,10.00){\line(-1,1){4.00}}
\put(96.00,44.00){\line(-1,-1){4.00}}
\put(62.00,10.00){\line(1,1){4.00}}
\put(80.00,65.00){\makebox(0,0)[cb]{$\l_r$}}
\put(88.00,65.00){\makebox(0,0)[cb]{$\l_{r-1}$}}
\put(96.00,65.00){\makebox(0,0)[cb]{$\l_{r-2}$}}
\put(105.00,65.00){\makebox(0,0)[cb]{$\cdots$}}
\put(122.00,65.00){\makebox(0,0)[cb]{$\l_{r-k+1}$}}
\put(130.00,65.00){\makebox(0,0)[cb]{$\l_{r-k}$}}
\put(73.00,65.00){\makebox(0,0)[cb]{$\l_{r+1}$}}
\put(30.00,65.00){\makebox(0,0)[cb]{$\l_{r+k}$}}
\put(38.00,65.00){\makebox(0,0)[cb]{$\l_{r+k-1}$}}
\put(59.00,65.00){\makebox(0,0)[cb]{$\cdots$}}
\put(25.00,60.00){\makebox(0,0)[rc]{$\l_k$}}
\put(25.00,52.00){\makebox(0,0)[rc]{$\l_{k-1}$}}
\put(25.00,44.00){\makebox(0,0)[rc]{$\l_{k-2}$}}
\put(25.00,36.00){\makebox(0,0)[rc]{$\vdots$}}
\put(25.00,10.00){\makebox(0,0)[rc]{$0$}}
\put(25.00,18.00){\makebox(0,0)[rc]{$\l_1$}}
\put(80.00,18.00){\makebox(0,0)[cc]{$\cdots$}}
\put(80.00,57.00){\makebox(0,0)[ct]{$+$}}
\put(72.00,49.00){\makebox(0,0)[ct]{$-$}}
\put(88.00,49.00){\makebox(0,0)[ct]{$+$}}
\put(96.00,41.00){\makebox(0,0)[ct]{$+$}}
\put(80.00,41.00){\makebox(0,0)[ct]{$-$}}
\put(122.00,15.00){\makebox(0,0)[ct]{$+$}}
\put(130.00,7.00){\makebox(0,0)[ct]{$+$}}
\put(114.00,7.00){\makebox(0,0)[ct]{$-$}}
\put(64.00,41.00){\makebox(0,0)[ct]{$+$}}
\end{picture}
\caption{The extended twisted TPG for
$U_q(A_{2l}^{(2)})$ for the product
$V_0(\l_k)\otimes V_0(\l_r)$ $(k\leq r,~r+k\leq l)$ .
The nodes correspond to representations whose highest weight
is given by
the sum of the weight labeling the column and the weight labeling
the row. The $\pm$ indicate the parity. The parities are equal
along the northwest-southeast diagonals and they alternate along
the northeast-southwest diagonals.
\label{fig1}}
\end{figure}
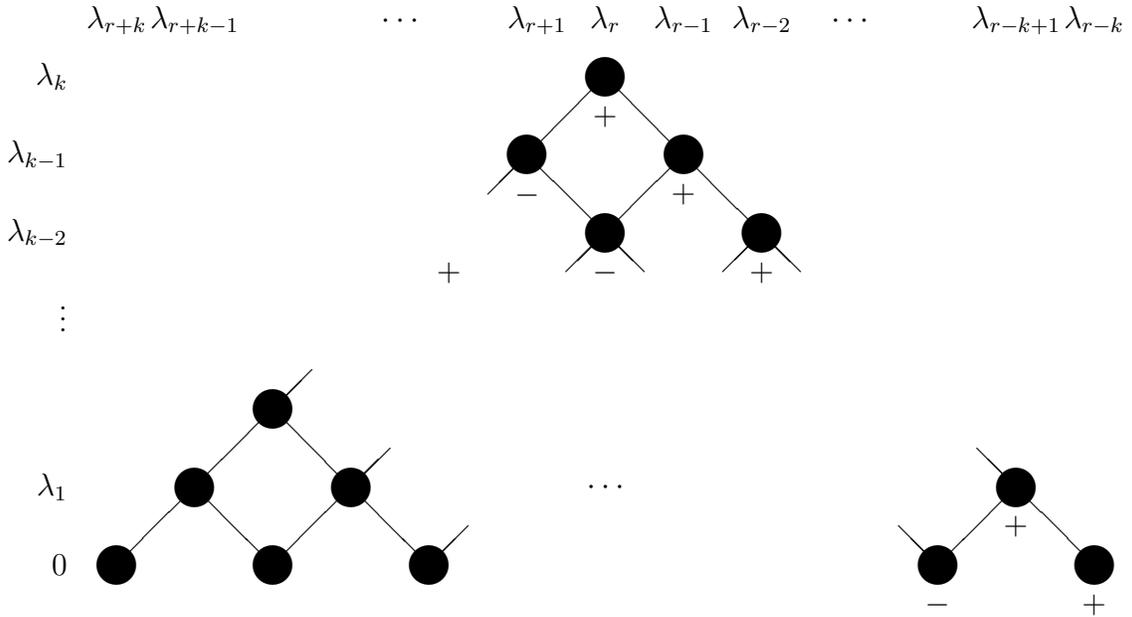

Now for $1\leq k\leq r\leq l$ we have the tensor product
decomposition
\[
V_0(\l_k)\otimes V_0(\l_r)=\bigoplus_{a=0}^k\bigoplus_{c=0}^a
V_0(\l_{c}+\l_{d})
\]
where
\[
d=\left\{\begin{array}{ll}
k+r-2a+c&\mbox{for }2a-c\geq r+k-l\\
2l+1-(k+r-2a+c)&\mbox{for }2a-c< r+k-l
\end{array}\right.
\]
which is multiplicity free (c.f. Appendix \ref{app2}).
The corresponding extended twisted tensor product graph is
consistent and quite different in topology to the the extended
untwisted TPG (which is inconsistent). We illustrate this
below with the case $r+k\leq l,~k\leq r$. In this case
$d=k+r-2a+c$ and we have the extended twisted TPG depicted in figure
\ref{fig1}.

To see that the extended twisted TPG in figure \ref{fig1} is
consistent, consider a typical closed loop:
\[
\unitlength=1.30mm
\linethickness{0.4pt}
\begin{picture}(89.00,40.00)(30,25)
\put(70.00,30.00){\circle*{4.00}}
\put(85.00,45.00){\circle*{4.00}}
\put(55.00,45.00){\circle*{4.00}}
\put(70.00,60.00){\circle*{4.00}}
\put(70.00,60.00){\line(1,-1){15.00}}
\put(85.00,45.00){\line(-1,-1){15.00}}
\put(70.00,30.00){\line(-1,1){15.00}}
\put(55.00,45.00){\line(1,1){15.00}}
\put(70.00,60.00){\line(0,0){0.00}}
\put(70.00,57.00){\makebox(0,0)[ct]{$+$}}
\put(82.00,45.00){\makebox(0,0)[rc]{$+$}}
\put(58.00,45.00){\makebox(0,0)[lc]{$-$}}
\put(70.00,33.00){\makebox(0,0)[cb]{$-$}}
\put(70.00,63.00){\makebox(0,0)[cb]{$\l_c+\l_d$}}
\put(89.00,45.00){\makebox(0,0)[lc]{$\l_{c-1}+\l_{d-1}$}}
\put(51.00,45.00){\makebox(0,0)[rc]{$\l_{c-1}+\l_{d+1}$}}
\put(70.00,27.00){\makebox(0,0)[ct]{$\l_{c-2}+\l_d$}}
\end{picture}
\]
where we have indicated the relative parities of the
vertices. Using the fact that on $V_0(\l_c+\l_d)$ the
universal Casimir element of $\lie_0$ takes the eigenvalue
\[
C_{c,d}=(c+d)(2l+2-c)-(d+1)(d-c),
\]
it is easily seen that
\[
&&
C_{c,d}-C_{c-1,d-1}=C_{c-1,d+1}-C_{c-2,d}=2(2l+3-c-d),
\non&&
C_{c,d}-C_{c-1,d+1}=C_{c-1,d-1}-C_{c-2,d}=2(d-c+2).
\]
This implies that the extended twisted TPG is consistent, i.e.
that the recursion relations \Eq{rec} give the
same result independent of the path along which one recurses.

We are now in a position to write down our solution to Jimbo's
equation and thus to the quantum Yang-Baxter equation arising from
the above extended twisted TPG:
\[
\check{R}^{\l_k,\l_r}(u)=
\sum_{a=0}^k\sum_{c=0}^a
\prod_{i=a}^{k-1}\<k+r-2i\>_-
\prod_{j=1}^{a-c}\<n-r-k+2j\>_+
\,\tp^{\l_k,\l_r}_{\l_c+\l_{k+r-2a+c}}
\]


\subsection{R-matrices for $U_q(A_{2l-1}^{(2)})$
\label{sectjimbo2}}

This is the case of the second line in table \ref{tab1}, i.e.
$\lie=A_{2l-1}=sl(2l)$, $\lie_0=C_l=
sp(2l)$ and $\theta_0=\e_1+\e_2=\l_2$.

Starting with the vector irrep $V_0(\l_1)$ of $U_q(\lie_0)$
(and also of $U_q(\lie)$) we have the following twisted
TPG for $V(\l_1)\otimes V(\l_1)$
\[\label{ttpg3}
\unitlength=1mm
\linethickness{0.4pt}
\begin{picture}(102.60,12.60)(40,13)
\put(60.00,20.00){\circle*{5.20}}
\put(80.00,20.00){\circle*{5.20}}
\put(100.00,20.00){\circle*{5.20}}
\put(60.00,20.00){\line(1,0){40.00}}
\put(60.00,16.00){\makebox(0,0)[ct]{$0$}}
\put(80.00,16.00){\makebox(0,0)[ct]{$\l_2$}}
\put(100.00,16.00){\makebox(0,0)[ct]{$2\l_1$}}
\put(60.00,24.00){\makebox(0,0)[cb]{$-$}}
\put(80.00,24.00){\makebox(0,0)[cb]{$-$}}
\put(100.00,24.00){\makebox(0,0)[cb]{$+$}}
\end{picture}
\]
which has quite a different topology to the untwisted TPG
\[\label{tpg3}
\unitlength=1mm
\linethickness{0.4pt}
\begin{picture}(102.60,12.60)(40,13)
\put(60.00,20.00){\circle*{5.20}}
\put(80.00,20.00){\circle*{5.20}}
\put(100.00,20.00){\circle*{5.20}}
\put(60.00,20.00){\line(1,0){40.00}}
\put(60.00,16.00){\makebox(0,0)[ct]{$0$}}
\put(80.00,16.00){\makebox(0,0)[ct]{$2\l_1$}}
\put(100.00,16.00){\makebox(0,0)[ct]{$\l_2$}}
\put(60.00,24.00){\makebox(0,0)[cb]{$-$}}
\put(80.00,24.00){\makebox(0,0)[cb]{$+$}}
\put(100.00,24.00){\makebox(0,0)[cb]{$-$}}
\end{picture}
\]
Because $V_0(2\l_1)$ appears as an extremal node on the twisted
TPG \Eq{ttpg3}, it is affinizable. Continuing in this way it
is easily seen that $V_0(a\l_1)$ is affinizable for any positive
integer $a$. We have the following $U_q(\lie_0)$-module decomposition
of the tensor product of any two such representations
\[
V_0(k\l_1)\otimes V_0(r\l_1)=\bigoplus_{a=0}^k\bigoplus_{b =0}^a
V_0((k+r-2a)\l_1+b\l_2),~~~k\leq r.
\]
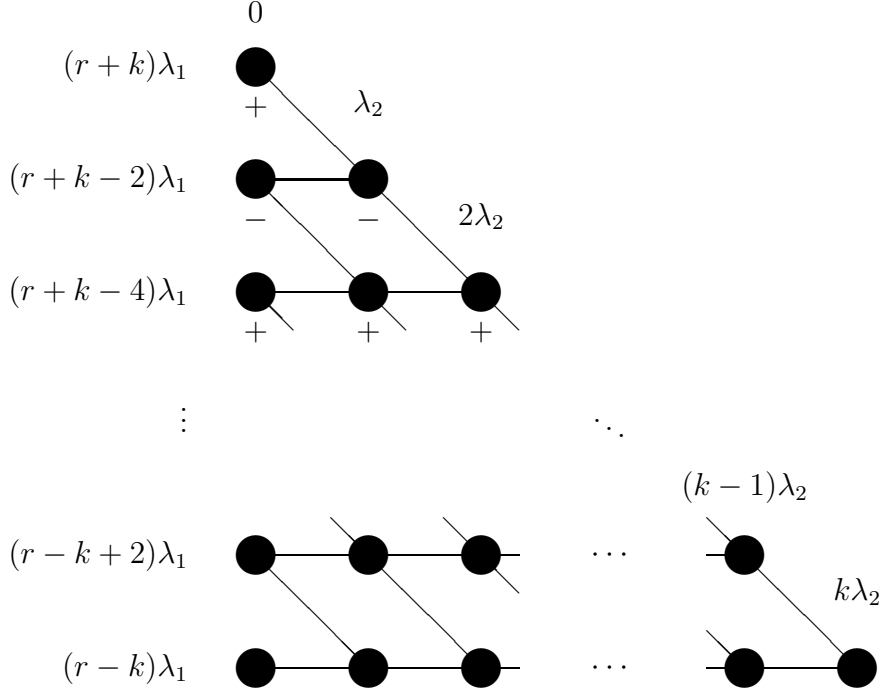
\begin{figure}[ht]
\unitlength=1.00mm
\linethickness{0.4pt}
\begin{picture}(129.60,92.00)
\put(47.00,86.00){\line(1,-1){35.00}}
\put(77.00,56.00){\line(-1,0){30.00}}
\put(67.00,51.00){\line(-1,1){20.00}}
\put(47.00,71.00){\line(1,0){15.00}}
\put(47.00,56.00){\line(1,-1){5.00}}
\put(47.00,6.00){\line(1,0){35.00}}
\put(47.00,21.00){\line(1,0){35.00}}
\put(47.00,21.00){\line(1,-1){15.00}}
\put(77.00,6.00){\line(-1,1){20.00}}
\put(127.00,6.00){\line(-1,0){20.00}}
\put(112.00,6.00){\line(-1,1){5.00}}
\put(127.00,6.00){\line(-1,1){20.00}}
\put(112.00,21.00){\line(-1,0){5.00}}
\put(72.00,26.00){\line(1,-1){10.00}}
\put(47.00,86.00){\circle*{5.20}}
\put(47.00,71.00){\circle*{5.20}}
\put(62.00,71.00){\circle*{5.20}}
\put(62.00,56.00){\circle*{5.20}}
\put(47.00,56.00){\circle*{5.20}}
\put(77.00,56.00){\circle*{5.20}}
\put(77.00,21.00){\circle*{5.20}}
\put(62.00,21.00){\circle*{5.20}}
\put(47.00,21.00){\circle*{5.20}}
\put(47.00,6.00){\circle*{5.20}}
\put(62.00,6.00){\circle*{5.20}}
\put(77.00,6.00){\circle*{5.20}}
\put(112.00,6.00){\circle*{5.20}}
\put(127.00,6.00){\circle*{5.20}}
\put(112.00,21.00){\circle*{5.20}}
\put(47.00,82.00){\makebox(0,0)[ct]{$+$}}
\put(47.00,67.00){\makebox(0,0)[ct]{$-$}}
\put(47.00,52.00){\makebox(0,0)[ct]{$+$}}
\put(62.00,67.00){\makebox(0,0)[ct]{$-$}}
\put(62.00,52.00){\makebox(0,0)[ct]{$+$}}
\put(77.00,52.00){\makebox(0,0)[ct]{$+$}}
\put(94.00,21.00){\makebox(0,0)[cc]{$\cdots$}}
\put(94.00,6.00){\makebox(0,0)[cc]{$\cdots$}}
\put(38.00,86.00){\makebox(0,0)[rc]{$(r+k)\l_1$}}
\put(38.00,71.00){\makebox(0,0)[rc]{$(r+k-2)\l_1$}}
\put(38.00,56.00){\makebox(0,0)[rc]{$(r+k-4)\l_1$}}
\put(38.00,21.00){\makebox(0,0)[rc]{$(r-k+2)\l_1$}}
\put(38.00,6.00){\makebox(0,0)[rc]{$(r-k)\l_1$}}
\put(62.00,80.00){\makebox(0,0)[cb]{$\l_2$}}
\put(77.00,65.00){\makebox(0,0)[cb]{$2\l_2$}}
\put(112.00,29.00){\makebox(0,0)[cb]{$(k-1)\l_2$}}
\put(127.00,15.00){\makebox(0,0)[cb]{$k\l_2$}}
\put(38.00,40.00){\makebox(0,0)[rc]{$\vdots$}}
\put(94.00,39.00){\makebox(0,0)[cc]{$\ddots$}}
\put(47.00,92.00){\makebox(0,0)[cb]{$0$}}
\end{picture}
\caption{The extended twisted TPG for $U_q(A_{2l-1}^{(2)})$
for the tensor product $V(k\l_1)
\otimes V(r\l_1)$.
The nodes correspond to modules whose highest weight is the sum of the
weight labeling the column and the weight labeling the row. Some of
the parities are indicated below the nodes. The parities are the same
along the horizontal and they alternate along the vertical.
\label{fig2}}
\end{figure}

The corresponding extended twisted TPG is shown in figure \ref{fig2}.
To see that this graph is consistent we have to consider the
closed loops of the form
\[
\unitlength=1.30mm
\linethickness{0.4pt}
\begin{picture}(93.00,28.00)(30,20)
\put(60.00,40.00){\circle*{4.00}}
\put(75.00,40.00){\circle*{4.00}}
\put(75.00,25.00){\circle*{4.00}}
\put(89.00,25.00){\circle*{4.00}}
\put(60.00,40.00){\line(1,0){15.00}}
\put(75.00,40.00){\line(1,-1){14.00}}
\put(89.00,26.00){\line(-1,0){14.00}}
\put(75.00,26.00){\line(-1,1){15.00}}
\put(56.00,40.00){\makebox(0,0)[rc]{$c\l_1+(d-1)\l_2$}}
\put(79.00,40.00){\makebox(0,0)[lc]{$c\l_1+d\l_2$}}
\put(93.00,25.00){\makebox(0,0)[lc]{$(c-2)\l_1+(d+1)\l_2$}}
\put(71.00,25.00){\makebox(0,0)[rc]{$(c-2)\l_1+d\l_2$}}
\put(60.00,37.00){\makebox(0,0)[ct]{$+$}}
\put(75.00,37.00){\makebox(0,0)[ct]{$+$}}
\put(75.00,22.00){\makebox(0,0)[ct]{$-$}}
\put(89.00,22.00){\makebox(0,0)[ct]{$-$}}
\end{picture}
\]
where we have indicated the relative parities of the vertices.
Using the fact that on $V_0(c\l_1+d\l_2)$ the quadratic Casimir
element of $\lie_0$ takes the eigenvalue
\[
C_{c,d}=(c+d)(n+c+d)+(n-2+d)d,
\]
it is easily seen that the above loop is consistent:
\[
&&
C_{c,d}-C_{c,d-1}=C_{c-2,d+1}-C_{c-2,d}=2(n-2+c+2d),
\non&&
C_{c,d}-C_{c-2,d+1}=C_{c,d-1}-C_{c-2,d}=2c.
\]

We can now read off the R-matrix from the extended twisted TPG
\[
\check{R}^{k\l_1,r\l_1}(u)&=&
\sum_{a=0}^k\sum_{b=0}^a
\prod_{i=1}^{(a-b)}\<n+k+r-2i\>_+
\prod_{j=1}^a\<k+r+2-2j\>_-
\non&&~~\times
\tp^{k\l_1,r\l_1}_{(k+r-2a)\l_1+b\l_2}.
\]
Again it can be seen that the above eigenvalues
of $\check{R}(u)$ are quite different to those arising from the
untwisted case and thus give rise to new solutions to the
quantum Yang-Baxter equation.


\subsection{R-matrices for $U_q(D_{l+1}^{(2)})$
\lab{sectjimbo3}}

This is the case of $\lie=D_{l+1}=so(2l+2),~\lie_0=B_l=so(2l+1)$
and $\theta_0=\e_1=\l_1$. We set $n=2l+1$.

The fundamental spinor irrep $V_0(\l_l),~\l_l=\half(\e_1+\cdots +
\e_l)$ of $U_q(\lie_0)$ is undeformed and affinizable and also carries
the spinor irreps $V(\l_l)$ and $V(\l_{l+1})$ of $U_q(\lie)$.
We have the tensor product decomposition, for $a\in\ZZ_+$,
\[
V_0(\l_l)\otimes V_0(a\l_l)=
\bigoplus_{k=0}^{l-1} V_0\left(\l_k+(a-1)\l_l\right)
\oplus V_0\left((a+1)\l_l\right)
\]
with the corresponding twisted TPG
\[\label{ttpg4}
\unitlength=1mm
\linethickness{0.4pt}
\begin{picture}(112.00,14.00)(4,10)
\put(20.00,18.00){\circle*{4.00}}
\put(35.00,18.00){\circle*{4.00}}
\put(65.00,18.00){\circle*{4.00}}
\put(80.00,18.00){\circle*{4.00}}
\put(110.00,18.00){\circle*{4.00}}
\put(50.00,18.00){\circle*{4.00}}
\put(20.00,18.00){\line(1,0){67.00}}
\put(104.00,18.00){\line(1,0){6.00}}
\put(96.00,18.00){\makebox(0,0)[cc]{$\cdots$}}
\put(20.00,14.00){\makebox(0,0)[ct]
{\scriptsize$(a+1)\l_l$}}
\put(35.00,14.00){\makebox(0,0)[ct]
{\scriptsize$\mu_1$}}
\put(50.00,14.00){\makebox(0,0)[ct]
{\scriptsize$\mu_2$}}
\put(65.00,14.00){\makebox(0,0)[ct]{$\cdots$}}
\put(110.00,14.00){\makebox(0,0)[ct]
{\scriptsize$(a-1)\l_l$}}
\put(20.00,21.00){\makebox(0,0)[cb]{$+$}}
\put(35.00,21.00){\makebox(0,0)[cb]{$-$}}
\put(50.00,21.00){\makebox(0,0)[cb]{$-$}}
\put(65.00,21.00){\makebox(0,0)[cb]{$+$}}
\put(80.00,21.00){\makebox(0,0)[cb]{$+$}}
\end{picture}
\]
where we have indicated the parity of the vertices
and introduced $\mu_i=(a-1)\l_l+\l_{l-i}$.
Note that his graph has a very different topology from the
corresponding untwisted TPG given in Fig. 1 of \ct{Del94a}. Because
$(a+1)\l_l$ is on an extremal node of this twisted TPG it
then follows by induction that $V_0(a\l_l)$ is affinizable
for any positive integer $a$. The
eigenvalue of the universal Casimir element of $\lie_0$
on $V_0(\l_k+(a-1)\l_l)$ is given by
\[
C\left(\l_k+(a-1)\l_l\right)&=&k(n+a-k-1)+\frac{1}{4}l(a-1)(n+a-2),~~
1\leq k\leq l-1,\non
C\left((a+1)\l_l\right)&=&l(n+a-l-1)+\frac{1}{4}l(a-1)(n+a-2)
\]
Using this we can read off the R-matrix from the twisted
TPG \Eq{ttpg4}:
\[
\check{R}^{\l_l,a\l_l}(u)=
\sum_{k=0}^{l-1}\rho_k(u)\,\tp^{\l_l,a\l_l}_{\l_k+(a-1)\l_l}
+\tp^{\l_l,a\l_l}_{(a+1)\l_l}
\]
with the eigenvalues $\rho_k(u)$ given by
\[
\rho_k(u)=\prod_{i=1}^{l-k}\left\langle\half(a-1)+i
\right\rangle_{(-1)^i}.
\]

We now proceed to the general case $V_0(a\l_l)\otimes V_0(b\l_l),~
0\leq a\leq b\in\ZZ$. In view of appendix \ref{app2} we now have
the multiplicity-free tensor product decomposition
\[
V_0(a\l_l)\otimes V_0(b\l_l)=\bigoplus_\L V_0(\L+(b-a)\l_l),
\]
where the sum is over all dominant weights
$\L=(\L_1,\L_2,\cdots,\L_l)$ satisfying
\[
a\geq\L_1\geq\L_2\geq\cdots\geq\L_l\geq 0,~~\L_i\in\ZZ,
\]
each such weight appearing exactly once.

The extended twisted TPG will typically be
$l$-dimensional and we can therefore not draw a diagramm for
it. However to determine wether it is consistent, it is
sufficient to look at closed loops
of the form (labeling the vertices by the weight $\L$ since
$a$ and $b$ are here fixed and thus redundant labels):
\[\label{loop}
\unitlength=1.30mm
\linethickness{0.4pt}
\begin{picture}(89.00,40.00)(30,25)
\put(70.00,30.00){\circle*{4.00}}
\put(85.00,45.00){\circle*{4.00}}
\put(55.00,45.00){\circle*{4.00}}
\put(70.00,60.00){\circle*{4.00}}
\put(70.00,60.00){\line(1,-1){15.00}}
\put(85.00,45.00){\line(-1,-1){15.00}}
\put(70.00,30.00){\line(-1,1){15.00}}
\put(55.00,45.00){\line(1,1){15.00}}
\put(70.00,60.00){\line(0,0){0.00}}
\put(70.00,63.00){\makebox(0,0)[cb]{$\L$}}
\put(89.00,45.00){\makebox(0,0)[lc]{$\L-\e_i$}}
\put(51.00,45.00){\makebox(0,0)[rc]{$\L-\e_j$}}
\put(70.00,27.00){\makebox(0,0)[ct]{$\L-\e_i-\e_j$}}
\end{picture}
\]
To show that all the edges in this loop really exist,
i.e. that \Eq{edge} is satisfied, one uses
the following theorem proven in \ct{Gou89}.
\begin{theorem}
Suppose $\l,\mu$ are dominant weights and $\nu$ is a weight
Weyl group conjugate to $\l$. If $\mu+\nu$ is dominant then
$V(\mu+\nu)$ occurs exactly once in $V(\l)\otimes V(\mu)$.
\end{theorem}

To obtain the relative (to the top vertex) parities of the
vertices of the closed loops
it is necessary to determine which irreps of $so(2l+1)$
belong to the same irrep of $so(2l+2)$ (which will all have the
{\em same} parity). As seen in appendix \ref{app2}, two such
irreps $V_0(\L+(b-a)\l_l)$ and $V_0(\L'+(b-a)\l_l)$ belong to the
same irrep of $so(2l+2)$ iff $\L_i=\L_i'$ for $l+i$ odd. From this
it follows that the difference in parity along any edge in
\Eq{loop} is equal to the difference in parity along the opposite
edge.
Using the eigenvalue formula for the universal Casimir element
of $\lie_0$ on $V_0(\L+(b-a)\l_l)$
\[
C_\L=\sum_{i=1}^l\L_i(\L_i+b-a+n-2i)+\frac{1}{4}l(b-a)(b-a+n-1),
\]
we see that also the difference of the Casimirs are equal along
opposite edges:
\[
C_\L-C_{\L-\e_i}=C_{\L-\e_j}-C_{\L-\e_i-\e_j}=
2(\L_i-i)+b-a+n-1,~~i\neq j.
\]
{}From these facts it follows that the extended twisted TPG is
consistent.

We are now in a position to determine the eigenvalues $\rho_\L(u)$
in the expression for the R-matrix
\[
\check{R}^{a\l_l,b\l_l}(u)=
\sum_\L\,\rho_{\L}(u)\,
\check{P}^{a\l_l,b\l_l}_{\L+(b-a)\l_l}.
\]
We start with the first component $\L_1$ and proceed along the following
path:
\[\label{path}
\unitlength=1mm
\linethickness{0.4pt}
\begin{picture}(92.00,12.00)(20,12)
\put(30.00,20.00){\circle*{4.00}}
\put(45.00,20.00){\circle*{4.00}}
\put(75.00,20.00){\circle*{4.00}}
\put(90.00,20.00){\circle*{4.00}}
\put(30.00,20.00){\line(1,0){21.00}}
\put(69.00,20.00){\line(1,0){21.00}}
\put(60.00,20.00){\makebox(0,0)[cc]{$\cdots$}}
\put(30.00,17.00){\makebox(0,0)[ct]{$\L$}}
\put(45.00,17.00){\makebox(0,0)[ct]{$\L+\e_1$}}
\put(90.00,17.00){\makebox(0,0)[lt]{$\L+(a-\L_1)\e_1$}}
\end{picture}
\]
We then proceed in this way component by component. By this means we
arrive at the formula
\[
\rho_\L(u)=\prod_{i=1}^l\prod_{k=\L_i}^{a-1}
\left\langle k-i+l+1+\half(b-a)\right\rangle_{(-1)^{l+i}}
\]
If $\L_i=a$ for some $i$ then the $i$th term is understood to
contribute 1 to the product. It is easily seen that this
formula reduces to our previous one when $a=b=1$.

\vspace{5mm}

\noindent{\bf Acknowledgements}
G.W.D. thanks the Deutsche Forschungsgemeinschaft for a
Habilitationsstipendium. Y.Z.Z. is financially supported by the
Kyoto University Foundation.


\appendix

\section{Kac generators\label{app1}}

In this appendix we give detailed expressions for the
generators $E_i,F_i$ and $H_i$ of \Eq{lierels}. We will
write them in terms of the
familiar basis elements $e_{ij},~i\leq i,j\leq n$ of  $gl(n)$
which satisfy the Lie bracket
\[
[e_{ij},e_{kl}]=\delta_{kj}e_{il}-\delta_{il}e_{kj}.
\]
As our invariant bilinear form on $gl(n)$ we take
\[\label{glnform}
(e_{ij},e_{kl})=\half\d_{jk}\d_{il}.
\]
We choose a basis $\{\e_i|1\leq i\leq n\}$ for the root space,
i.e. the dual space
to the Cartan subalgebra of $gl(n)$ such that $\sqrt{2}e_{ii}$ is
paired with $\e_i$. Then the bilinear form \Eq{glnform}
induces the scalar product $(\e_i,\e_j)=\d_{ij}$.

The diagram automorphism $\s$ of order $k=2$ by which we will
twist is given by
\[
\s(e_{ij})=(-1)^{i+j+1}e_{\bar{\jmath}\bar{\imath}},~~~~
\bar{\imath}=n+1-i.
\]
The fixed point subalgebra $\lie_0$
is generated by the linear combinations
\[
a_{ij}=e_{ij}-(-1)^{i+j}e_{\bar{\jmath}\bar{\imath}}.
\]
Its Cartan subalgebra $H$ is spanned by the $a_{ii},~1\leq i\leq
\left[\frac{n}{2}\right]$.
The subspace $\lie_1$ is spanned by the elements
\[
b_{ij}=e_{ij}+(-1)^{i+j}e_{\bar{\jmath}\bar{\imath}}-
\frac{2}{n}\d_{ij}I_1,~~~I_1=\sum_{i=1}^n e_{ii}.
\]

\subsection{$\lie=A_{2l}=sl(2l+1),~\lie_0=B_l=so(2l+1)$
\label{app11}}

In this case we perform the above construction with $n=2l+1$.
Our simple generators are given by
\[
&&E_i=a_{i,i+1},~~F_i=a_{i+1,i},~~H_i=a_{i,i}-a_{i+1,i+1},~~~
1\leq i< l,
\non&&
E_l=a_{l,l+1},~~F_l=a_{l+1,l},~~H_l=a_{l,l},
\non&&
E_0=\frac{1}{\sqrt{2}} b_{n1}=\sqrt{2}a_{2l+1,1},~~
F_0=\sqrt{2}a_{1,2l+1},~~
H_0=-2a_{11},
\]
with the corresponding simple roots
\[
\rootb_i=\e_i-\e_{i+1},~~(1\leq i<l),~~
\rootb_l=\e_l,~~\rootb_0=-2\e_1.
\]
They can be checked to satisfy the defining relations \Eq{lierels}.
The invariant bilinear form \Eq{glnform} induces the form \Eq{form}.

\noindent{\bf Note:} Our notation here differs from that of Kac \ct{Kac90} who
interchanges indices $i=0$ and $i=l$. We prefer our notation since
it conforms to the usual notation in the untwisted case.

\subsection{$\lie=A_{2l-1}=sl(2l),~\lie_0=C_l=sp(2l)$
\label{app12}}

In this case $n=2l$ and our simple generators are
\[
&&E_i=a_{i,i+1},~~F_i=a_{i+1,i},~~H_i=a_{i,i}-a_{i+1,i+1},~~~
1\leq i< l,
\non&&
E_l=\frac{1}{\sqrt{2}}a_{l,l+1},~~
f_l=\frac{1}{\sqrt{2}}a_{l+1,l},~~
H_l=2a_{l,l},
\non&&
E_0=b_{2l-1,1},~~
F_0=b_{1,2l-1},~~
H_0=-(a_{11}+a_{22}),
\]
with the corresponding simple roots
\[
\rootb_i=\e_i-\e_{i+1},~~(1\leq i<l),~~
\rootb_l=2\e_l,~~\rootb_0=-(\e_1+\e_2).
\]
They can be checked to satisfy the defining relations \Eq{lierels}.
The invariant bilinear form \Eq{glnform} induces the form \Eq{form}.

\subsection{$\lie=D_{l+1}=so(2l+2),~\lie_0=B_l=so(2l+1)$
\label{app13}}

In this case we embedd $so(2l+1)$ and $so(2l+2)$ into $gl(2l+2)$.
The $\lie_0=so(2l+1)$ generators are
\[
a_{i,j}=e_{ij}-(-1)^{i+j}e_{\bar{\jmath}\bar{\imath}},~~~~
1\leq i,j\leq 2l+1,
\]
where $\bar{\imath}=2l+2-i$.
The extra generators which span $\lie_1$ and complete $so(2l+1)$
to $so(2l+2)$ are
\[
a_{i,2l+2}=e_{i,2l+2}-(-1)^i e_{2l+2,\bar{\imath}},~~~~
1\leq i\leq 2l+1
\]
As the set of simple generators for $\lie=so(2l+2)$ we take
\[
&&E_i=a_{i,i+1},~~F_i=a_{i+1,i},~~H_i=a_{i,i}-a_{i+1,i+1},~~~
1\leq i< l,
\non&&
E_l=a_{l,l+1},~~
f_l=a_{l+1,l},~~
H_l=a_{l,l},
\non&&
E_0=b_{2l+1,2l+2},~~
F_0=b_{2l+2,2l+1},~~
H_0=-a_{11},
\]
with the corresponding simple roots
\[
\rootb_i=\e_i-\e_{i+1},~~(1\leq i<l),~~
\rootb_l=\e_l,~~\rootb_0=-\e_1.
\]
They can be checked to satisfy the defining relations \Eq{lierels}.
The invariant bilinear form \Eq{glnform} induces the form \Eq{form}.


\section{Decompositions and Branching rules
\label{app2}}

If $\pi$ is any finite dimensional irrep of $\lie$ on a space
$V$ then we obtain
a representation $\hat{\pi}$ of $\tlie$ on the loop space $V(t)$ with
the prescription
\[
&&
\hat{\pi}(e_i)=\pi(E_i),~~~
\hat{\pi}(f_i)=\pi(F_i),~~~
\hat{\pi}(h_i)=\pi(H_i),~~~i\leq i\leq l,
\non&&
\hat{\pi}(e_0)=t \pi(E_0),~~~
\hat{\pi}(f_0)=t^{-1}\pi(F_0),~~~
\hat{\pi}(h_0)=\pi(H_0).
\]
Here we wish to consider the
\minimal irreps of $\lie$ which, by definition, are also irreps of
$\lie_0$. From the known $\lie\downarrow\lie_0$ branching rules
this restricts us to the irreps of $\lie$ with highest weights
given in table \ref{tab2}. These \minimal irreps are
of interest because they can always be quantized to give rise to
irreps of the twisted quantum affine algebra $U_q(\tlie)$, as we
have seen in the paper. Our aim here is to show that in the
decomposition of the tensor product of any two \minimal irreps
all irreps of $\lie_0$ occur with at most unit multiplicity. At
the same time we shall deduce the tensor product branching rules
used in the paper. As before it is convenient to treat each of
the three families separately.


\subsection{$A_{2l}^{(2)}$
\label{app21}}

This is the case with $\lie=A_{2l}=sl(n=2l+1),~\lie_0=B_l=so(2l+1)$.
For $k\leq l$ we have seen that $V_0(\l_k)=V(\l_k)$,
where $\l_k=\sum_{i=1}^k\e_i$, is \minimal.
We first consider the decomposition of the
tensor product of two such irreps into irreps of $\lie$
\[
V(\l_k)\otimes V(\l_r)=\bigoplus_{a=0}^k V(\l_a+\l_b),~~~
b=k+r-a,~~1\leq k\leq r\leq l.
\]
We have
\[
\mbox{dim} V(\l_a+\l_b)=\frac{b-a+1}{n+1}
\binco{n+1}{a}\binco{n+1}{b+1}.
\]

{}From well established techniques (e.g. Young diagrams \ct{Ham}
or the quasi-spin formalism \ct{Jud71,Bie81,Hec73}) we deduce
the $\lie\downarrow\lie_0$ branching rules
(here and below $a\wedge b=\mbox{min}(a,b)$)
\[
V(\l_a+\l_b)=\bigoplus_{c =0}^{a\wedge (n-b)}
V_0(\l_{c }+\l_{d }),~~~
d =(b-a+c )\wedge(n-(b-a+c )).
\]
This decomposition is consistent with the dimension
formula
\[
\mbox{dim}V_0(\l_{c }+\l_{d })=
\frac{(1+d -c )(n+1-c -d )}{(n+1)(n+2)}
\binco{n+2}{c }\binco{n+2}{d +1}.
\]
For the case at hand $1\leq a\leq k\leq r\leq l$ from which it
follows that $a\wedge (n-b)=a$ and we obtain the $\lie_0$-module
decomposition
\[
V_0(\l_k)\otimes V_0(\l_r)&=&
\bigoplus_{a=0}^k\bigoplus_{c =0}^a V_0(\l_{c }+\l_{d }),
\non
d &=&(k+r-2a+c )\wedge(n-(k+r-2a+c )).
\]
To see that this decomposition is indeed multiplicity free
suppose that the $\lie_0$-module $V_0(\l_c +\l_d )$ occured twice.
This could only happen if there existed $a,a'$ such that
$(k+r-2a+c )=n-(k+r-2a'+c )$ which is however impossible
because $n=2l+1$ is odd.


\subsection{$A_{2l-1}^{(2)}$
\label{app22}}

This is the case of $\lie=A_{2l-1}=sl(n=2l),~\lie_0=C_l=sp(2l)$.
In this case $V(k\l_1)=V_0(k\l_1)$ is minimal and irreducible
as both an $\lie$-module and an $\lie_0$-module. For $k\leq r$
we have the $\lie$-module decomposition
\[
V(k\l_1)\otimes V(r\l_1)=\bigoplus_{a=0}^k V(b\l_1+a\l_2),~~~
b=k+r-2a.
\]
The $\lie
\downarrow\lie_0$ branching rules for such modules is
\[
V(b\l_1+a\l_2)=\bigoplus_{c =0}^a V_0(b\l_1+c\l_2),
\]
which again follows from
Young diagrams or the quasi-spin formalism  \ct{Jud71,Bie81}, \ct{Hec73}.
The decomposition is consistent with the dimension formulae
\[
\mbox{dim}V(b\l_1+a\l_2)&=&\frac{b+1}{n-1}\binco{a+b+n-1}{n-2}
\binco{a+n-2}{n-2},
\non
\mbox{dim}V_0(d\l_1+c\l_2)&=&\frac{(1+d )(2c +d +n-1)}{(n-1)(n-2)}
\binco{c +d +n-2}{n-3}\binco{c +n-3}{n-3}.
\]
We thus arrive at the irreducible $\lie_0$-module decomposition
\[
V_0(k\l_1)\otimes V_0(r\l_1)=\bigoplus_{a=0}^k\bigoplus_{c =0}^a
V_0((k+r-2a)\l_1+c\l_2),~~~k\leq r,
\]
which is easily seen to be multiplicity free.


\subsection{$D_{l+1}^{(2)}$
\label{app23}}

This is the case of $\lie=D_{l+1}=so(n=2l+2),~\lie_0=B_l=so(2l+1)$.

We first recall the PRV theorem \ct{Par67}. Let $\lie$ be a finite
dimensional simple Lie algebra with simple roots
$\a_i,1\leq i\leq l$ and corresponding Chevalley generators
$e_i,f_i,h_i,1\leq i\leq l$. We denote the set of dominant weights
of $\lie$ by $D_+$ and for $\l\in D_+$ we let $\Pi(\l)$ denote
the set of distinct weights in the finite dimensional irreducible
module $V(\l)$.
\begin{theorem}
For $\l,\mu\in D_+$ we have the tensor product decomposition
\[
V(\l)\otimes V(\mu)=\bigoplus_{\nu\in\Pi(\l)}
\mbox{dim}(V_{\nu,\mu}(\l))\,V(\mu+\nu).
\]
Here
\[
V_{\nu,\mu}(\l)=\left\{v\in V_\nu(\l)|e_i^{\<\mu,\a_i\>+1}v=0,
1\leq i\leq l\right\}
\]
where $V_\nu(\l)$ is the weight space consisting of weight vectors
of $V(\l)$ of weight $\nu\in\Pi(\l)$ and
\[
\<\mu,\a_i\>=\frac{2(\mu,\a_i)}{(\a_i,\a_i)}.
\]
\end{theorem}

To determine the $so(2l+1)$ branching rules for
$V_0(a\l_l)\otimes V_0(b\l_l),~a\leq b,~\l_l=(\half,\dots\half)$,
we use the fact that $V_0(a\l_l)$ is the carrier space for
parafermistatistics of order $a$ \ct{Bra72}.
{}From this it is known that $V_0(a\l_l)$
decomposes into a direct sum of irreps of the $gl(l)$
subalgebra with highest weight of the form $\L-a\l_l$ where
$\L=(\L_1,\dots,\L_l)$ satisfies
\[\label{Ls}
a\geq\L_1\geq\L_2\geq\cdots\geq\L_l\geq 0,~~\L_i\in\ZZ,
\]
each occuring exactly once. We deduce from this, together
with the PRV theorem,
the tensor product decomposition
\[
V_0(a\l_l)\otimes V_0(b\l_l)=\bigoplus_\L V_0(\L+(b-a)\l_l),
\]
where the sum is over all $\L$ satisfying \Eq{Ls}.

To obtain the correct parity pattern for the extended twisted
tensor product graph we need to investigate the $so(2l+2)$
branching law. It is convenient to consider the so called
statistical operator $\tilde{q}$ which takes a constant value
on the irreducible $gl(l)$ submodule with highest weight
$\L-a\l_l$ (with $\L$ as in \Eq{Ls}) given by
\[
\tilde{q}=-\sum_{i=1}^l(-1)^i\L_i.
\]
Then $\e_{l+1}^\vee\equiv a/2-\tilde{q}$ is the additional
Cartan generator of $so(2l+2)$ so that, in this representation,
$V_0(a\l_l)$ gives rise to an irreducible $so(2l+2)$ module
$V(a\hat{\l}_s)$ with highest weight $\hat{\l}_s$,
\[
\hat{\l}_s=\left(\half,\cdots,\half,(-1)^{l+1}\half\right)~~~~~
(l+1 \mbox{ components}).
\]
Again applying the PRV theorem we deduce the following $so(2l+2)$
branching law:
\[
V_0(a\l_l)\otimes V_0(b\l_l)=V(a\hat{\l}_s)\otimes V(b\hat{\l}_s)
=\bigoplus_{\hat{\L}} V\left(\hat{\L}+(b-a)\hat{\l}_s\right),
\]
where the sum is over all $so(2l+2)$ weights $\hat{\L}\in D_+$ of
the form
\[
\hat{\L}=\left\{
\begin{array}{ll}
\left(\hat{\L}_1,\hat{\L}_1,\hat{\L}_2,\hat{\L}_2,\cdots,
\hat{\L}_{l/2},\hat{\L}_{l/2},0\right),
& l \mbox{ even}\\
\left(a,\hat{\L}_1,\hat{\L}_1,\hat{\L}_2,\hat{\L}_2,\cdots,
\hat{\L}_{(l-1)/2},\hat{\L}_{(l-1)/2},0\right),
& l \mbox{ odd}
\end{array}\right.
\]
with $a\geq \hat{\L}_i\geq\hat{\L}_{i+1},~\forall i$.

\noindent{\bf Note:}
The above may be uniquely characterized as those weights $\hat{\L}=
(\L,0)$, where $\L\in D_0^+$ are those $so(2l+1)$ weights of
the form \Eq{Ls}, for which the corresponding statistical quantum
number $\tilde{q}$ takes its minimal (resp. maximal) value
$\tilde{q}=0$ (resp. $\tilde{q}=a$) when $l$ is even (resp. odd).

It follows from this that an irreducible $so(2l+1)$ module
$V_0(\L+(b-a)\l_l)$ belongs to the irreducible $so(2l+2)$
module $V(\hat{\L}+(b-a)\hat{\l}_s)$ where
\[
\hat{\L}_i=\left\{
\begin{array}{ll}
\L_i, & i+l \mbox{ odd}\\
\L_{i-1},& i+l \mbox{ even}
\end{array}\right.
\]
where we have set $\L_0\equiv a$. This leads to the parity pattern
used in the paper.


\section{Note on Parities\lab{app3}}

Here we derive the following result:
\begin{quote}
Two vertices $\nu\neq\nu'$ connected by an edge in the
twisted TPG (i.e., for which \Eq{edge} is satified)
must have the {\em same} parity if $V_0(\nu)$ and $V_0(\nu')$
belong to the {\em same} irreducible $\lie$-module while they
must have {\em opposite} parities if they belong to different
irreducible $\lie$-modules.
\end{quote}
We used this in section \ref{sectjimbo}
to show that the twisted tensor product graph is contained in
the extended twisted tensor product graph. The following
considerations are a generalization of those for the untwisted
case contained in Appendix B of \ct{Del94a} to which we refer the
reader for the details.

It is sufficient to consider the case $q=1$.
We introduce the {\em twisted} parity operator
\[
\tilde{P}\equiv R^{\l\mu}(1)|_{q=1},
\]
which satisfies $\tilde{P}^2=1$. As in eqs. (B.12) and (B.13)
of \ct{Del94a} we obtain the  equations
(omitting $\pi_\l$ and $\pi_\mu$ below)
\[\label{pp1}
\tilde{P}P^{(0)\l\mu}_\nu\left (e_0\otimes 1
  \right )P^{(0)\l\mu}_{\nu'}
  &=&P^{(0)\l\mu}_\nu\left (1\otimes e_0
  \right )P^{(0)\l\mu}_{\nu'}\tilde{P},
\non
\tilde{P}P^{(0)\l\mu}_\nu\left (1\otimes e_0
  \right )P^{(0)\l\mu}_{\nu'}
&=&P^{(0)\l\mu}_\nu\left (e_0\otimes 1
  \right )P^{(0)\l\mu}_{\nu'}\tilde{P}\,,~~~~\nu\neq\nu',
\]
where $P^{(0)\l\mu}_\nu=P^{\l\mu}_\nu|_{q=1}$. We will now show
that $\tilde{P}$ coincides with the normal (untwisted) parity
operator $P$ on $\lie_0$ defined in our previous work.
{}From \Eq{pp1} we deduce, since $\tilde{P}$ and $P^{(0)\l\mu}_\nu$
commute with the diagonal action of $\lie_0$, that, for $\nu\neq\nu'$
\[
\tilde{P}P^{(0)\l\mu}_\nu\left (a\otimes 1
  \right )P^{(0)\l\mu}_{\nu'}
  &=&P^{(0)\l\mu}_\nu\left (1\otimes a
  \right )P^{(0)\l\mu}_{\nu'}\tilde{P}\nonumber\\
\tilde{P}P^{(0)\l\mu}_\nu\left (1\otimes a
  \right )P^{(0)\l\mu}_{\nu'}
  &=&P^{(0)\l\mu}_\nu\left (a\otimes 1
  \right )P^{(0)\l\mu}_{\nu'}\tilde{P},~~~\forall a\in\lie_1
\]
since $\lie_1$ is irreducible under the adjoint action of $\lie_0$.
It follows that, for all $a\in\lie_1$,
$1\otimes a-a\otimes 1$ reverses the twisted
parity while $1\otimes a+a\otimes 1$ preserves it.
Then also, for $a,b\in\lie_1$,
\[
[1\otimes a+a\otimes 1,1\otimes b-b\otimes 1]=
1\otimes [a,b]-[a,b]\otimes 1
\]
must reverse the twisted parity and since $[\lie_1,\lie_1]
=\lie_0$ this implies that
$1\otimes a-a\otimes 1$ reverses
the parity for any $a\in\lie_0$.
But this is the defining property of the normal
(untwisted) parity operator $P$ so that both $P$ and $\tilde{P}$
have the same eigenvalues on the irreducible $\lie_0$ modules
$V_0(\nu)\subseteq V_0(\l)\otimes V_0(\mu)$. Thus we must have
$\tilde{P}=P$.

It follows that $1\otimes a+a\otimes 1$ preserves the usual
parity for all $a\in\lie$ while $1\otimes a-a\otimes 1$
reverses it. In particular, since all irreducible $\lie_0$
modules contained in a given irreducible $\lie$ module are
connected to one another by repeated application of the
generators $\D(a)=a\otimes 1+1\otimes a,~a\in\lie_1$, we
deduce that they must all have the {\em same} parity. Thus
two vertices $\nu\neq\nu'$ connected by an edge in the
twisted tensor product graph must have the {\em same}
parity if $V_0(\nu)$ and $V_0(\nu')$ belong to the {\em same}
irreducible $\lie$ module while they must have
{\em opposite} parities if they belong to different irreducible
$\lie$ modules since then
\[
0\neq P^{(0)\l\mu}_\nu(e_0\otimes 1)P^{(0)\l\mu}_{\nu'}
=\half P^{(0)\l\mu}_\nu\left(e_0\otimes 1
-1\otimes e_0\right)P^{(0)\l\mu}_{\nu'}.
\]



\begin{thebibliography}{99}
\bibitem {Art94} S. Artz, L. Mezincescu, R.I. Nepomechie,
   {\it Spectrum of transfer matrix for $U_q(B_n)$-invariant
   $A^{(2)}_{2n}$ open spin chain}, hep-th/9409130.
\bibitem {Art95}  S. Artz, L. Mezincescu, R.I. Nepomechie,
   {\it Analytic Bethe Ansatz for $A^{(2)}_{2n-1},B^{(1)}_n,C^{(1)}_n,
   D^{(1)}_n$ quantum-algebra-invariant open spin chains},
   hep-th/9504085.
\bibitem {Bax82}        R.J. Baxter,     {\it Exactly Solved Models
  in Statistical Mechanics}, Academic Press, London {\bf } (1982) .
\bibitem {Baz85}  V.V. Bazhanov,  {\it Trigonometric Solutions of Triangle
  Equations and Classical Lie Algebras}, Phys. Lett. {\bf 159B} (1985) 321.
\bibitem {Bie81} L.C. Biedenharn, J.D. Louck, {\it Anglular Momentum in
  Quantum Physics, Theory and Application}, Addison-Wesley, Reading, MA,
(1981).
\bibitem {Bra72} A.J. Bracken, H.S. Green, {\it Algebraic Identities for
  Parafermi Statistics of Given Order},  Nuovo Cim. {\bf 9} (1972) 349.
\bibitem {Del94a}  G.W. Delius, M.D. Gould, Y.-Z. Zhang,   {\it On the
  construction of trigonometric solutions of the Yang-Baxter equation},
  Nucl. Phys.  {\bf B432} (1994) 377.
\bibitem {Del95a}     G.W. Delius, M.D. Gould, J.R. Links, Y.-Z. Zhang,
  {\it On Type-I Quantum Affine Superalgebras}, hep-th/9408008,
  Int. J. Mod. Phys. {\bf A10} (1995) (in press).
\bibitem {Del95b}    G.W. Delius, Y.-Z. Zhang, {\it Finite dimensional
  representations of quantum affine Lie algebras}, J. Phys. {\bf A:}
  Math. Gen. {\bf 28} (1995) 1915.
\bibitem {Del95c}        G.W. Delius, {\it Exact S-matrices with quantum affine
  symmetry}, Nucl. Phys. {\bf B} to appear, hep-th/9503079.
\bibitem {Dri85}        V.G. Drinfel'd,  {\it Hopf algebras and the quantum
   Yang-Baxter equation}, Sov. Math. Dokl. {\bf 32} (1985) 254.
\bibitem {Gou89} M.D. Gould, {\it Casimir invariants and infinitesimal
  characters for semi-simple Lie algebras}, Rep. Math. Phys.
  {\bf 27} (1989) 73.
\bibitem {Gou94} M.D. Gould, J. Paldus, J. Cizek, {\it Quasi-Spin and
  the Pseudo-Orthogonal Group in the Hubbard Model}, Int. J. Quant. Chem.
  {\bf 50} (1994) 207.
\bibitem {Grimm} U. Grimm,  J. Phys. {\bf A:} Math. Gen. {\bf 27} (1994) 5897.
\bibitem {Ham} M. Hamermesh, {\it Group Theory and Its Applications}
  Addison-Wesley, Reading, MA, (1962).
\bibitem {Hec73} K.T. Hecht, Annu. Rev. Nucl. Sci. {\bf 23} (1973) 123.
\bibitem{Ize81}
  A.G. Izergin and V.E. Korepin, Commun. Math. Phys. {\it 79} (1981) 303.
\bibitem {Jim85}     M. Jimbo,   {\it A q-Difference Analogue of U(g) and
  the Yang-Baxter Equation}, Lett. Math. Phys. {\bf 10} (1985) 63.
\bibitem {Jim86}    M. Jimbo,  {\it Quantum R-matrix for the Generalized
  Toda System}, Commun. Math. Phys. {\bf 102} (1986) 537.
\bibitem {Jim89}    M. Jimbo,   {\it Introduction to the Yang-Baxter
  Equation}, Int. J. Mod. Phys. {\bf A4} (1989) 3759.
\bibitem {Jud71} B.R. Judd, Adv. At. Mol. Phys. {\bf 7} (1971) 251.
\bibitem {Kac90}   V.G. Kac,   {\it Infinite Dimensional Lie Algebras},
  Cambridge University Press {\bf } (1990) .
\bibitem {Kor93}     V.E. Korepin, N.M. Bogoliubov, A.G. Izergin,
   {\it Quantum inverse scattering method and correlations functions},
   Cambridge University Press {\bf } (1993) .
\bibitem {Kul81}   P.P. Kulish, N.Yu. Reshetikhin, E.K. Sklyanin,
  {\it Yang-Baxter equation and representation theory: I},
  Lett. Math. Phys. {\bf 5} (1981) 393.
\bibitem {Kun94} A. Kuniba, J. Suzuki, {\it Functional Relations and
   Analytic Bethe Ansatz for Twisted Quantum Affine Algebras},
   hep-th/9408135.
\bibitem {Par67} K.R. Parthasarathy, R. Ranga-Rao, V.S. Varadarajan,
  Ann. Math. {\bf 85} (1967) 383.
\bibitem{Res87}
  N. Yu. Reshetikhin, Lett. Math. Phys. {\it 14} (1987) 235.
\bibitem {Skl79} E.K. Sklyanin, L.A. Takhtadzhyan, L.D. Fadeev,
  {\it Quantum inverse problem method I}, Theor. Math. Phys. {\bf 40} (1979)
  194.
\bibitem {Zam79}    A.B Zamolodchikov, Al.B. Zamolodchikov,  {\it Factorized
  S-matrices in Two Dimensions as the Exact Solutions of Certain Relativistic
  Quantum Field Theory Models}, Ann. of Phys. {\bf 120} (1979) 253.
\bibitem {Zha91}    R.B. Zhang, M.D. Gould, A.J. Bracken,    {\it From
  Representations of the Braid Group to Solutions of the Yang-Baxter
  Equation}, Nucl. Phys. {\bf B354} (1991) 625.
\end{thebibliography}
\end{document}